\RequirePackage{ifpdf}
\ifpdf 
\documentclass[pdftex]{sigma}
\else
\documentclass{sigma}
\fi

\begin{document}

\renewcommand{\thefootnote}{$\star$}

\renewcommand{\PaperNumber}{076}

\FirstPageHeading

\ShortArticleName{Liouville Theorem for Dunkl Polyharmonic Functions}

\ArticleName{Liouville Theorem for Dunkl Polyharmonic Functions\footnote{This paper is a contribution to the Special
Issue on Dunkl Operators and Related Topics. The full collection
is available at
\href{http://www.emis.de/journals/SIGMA/Dunkl_operators.html}{http://www.emis.de/journals/SIGMA/Dunkl\_{}operators.html}}}

\Author{Guangbin  REN~$^{\dag\ddag}$ and Liang LIU~$^\dag$}

\AuthorNameForHeading{G.B. Ren  and L. Liu}

\Address{$^\dag$~Department of Mathematics,
 University of Science and Technology of China,\\
\hphantom{$^\dag$}~Hefei, Anhui 230026, P.R. China}
 \EmailD{\href{mailto:rengb@ustc.edu.cn}{rengb@ustc.edu.cn}, \href{mailto:xiaweije@mail.ustc.edu.cn}{xiaweije@mail.ustc.edu.cn}}

\Address{$^\ddag$~Departamento de Matem\'atica, Universidade de
Aveiro, P-3810-193, Aveiro, Portugal}

\ArticleDates{Received July 03, 2008, in f\/inal form October 30,
2008; Published online November 06, 2008}

\Abstract{Assume that   $f$ is Dunkl polyharmonic in $\mathbb{R}^n$
(i.e.\ $(\Delta _h)^p f=0$ for some integer~$p$, where $\Delta _h$ is
the Dunkl Laplacian associated to a root system $R$ and to a
multiplicity function~$\kappa$, def\/ined on $R$ and invariant with
respect to the f\/inite Coxeter  group).
Necessary and successful condition that $f$ is a polynomial of degree $\le s$ for $s\ge 2p-2$
is proved.
As a~direct corollary,  a Dunkl harmonic function bounded
above or below  is
 constant.}

 \vspace{-1mm}


\vspace{-1mm}

\Classification{33C52; 31A30; 35C10}

\vspace{-3mm}

\section{Introduction} The classical Liouville theorem for harmonic functions states that a
 harmonic function on $\mathbb{R}^{n}$ must be a constant if it is
 bounded or nonnegative.
 Nicolesco \cite{N} extended the Liouville theorem  to polyharmonic
functions with the  Pizetti formula as a starting point (see also~\cite{G}). Kuran~\cite{K},  Armitage~\cite{A1}, and  Futamura, Kishi,
and Mizuta \cite{FKM}  proved   further extensions and showed that
if~$f$ is a polyharmonic function on $\mathbb{R}^{n}$ and the growth
of the positive part of~$f$ is suitably restricted, then~$f$ must be
a polynomial. Their starting point is the Almansi decomposition
theorem for polyharmonic functions. We also refer to~\cite{L1,L2}
for the extension of Liouville theorems for conformally invariant
fully nonlinear equations. Recently, Gallardo and Godefroy~\cite{GG}
showed that if $f$ is a bounded  Dunkl harmonic function in
$\mathbb{R}^{n}$, then it is a constant. However, their approach
 is not
adaptable to Dunkl polyharmonic functions.

 The purpose of this article is to establish the Liouville theorem for   Dunkl
polyharmonic functions. To achieve this, we shall resort to  the
Almansi decomposition for Dunkl polyharmonic functions \cite{R}. As
a direct corollary
   of  our  results,
  a Dunkl harmonic function  bounded below or above is
    actually constant, which extends the corresponding result of  Gallardo and Godefroy  for the bounded
    case~\cite{GG}.
In the Dunkl ananlysis,  the multiplicity function  is usually
restricted to be non-negative. We shall discuss in the f\/inal section
the possible extension of our main result   to the case when the
multiplicity function is negative.

\section{Dunkl polyharmonic functions}

 For a
nonzero vector $v=(v_1,\ldots, v_n)\in\mathbb R^n$, the ref\/lection
$\sigma_{v}$ with respect to  the hyperplane orthogonal to $ v$  is
def\/ined by{\samepage
\[
\sigma_{v}{x}:
={x}-2\frac{\langle {x},{v} \rangle}{|v|^2}  v,\qquad
x=(x_1,\ldots,x_m)\in \mathbb R^m,
\]
where the symbol $\langle
 x,  y \rangle$ denotes the usual Euclidean
inner product and $| x|^2=\langle  x, x\rangle$.}

A root system $R$ is a f\/inite  set of nonzero vectors in $\mathbb
R^m$ such that $\sigma_{ v} R=R$ and $R\cap \mathbb R{ v}=\{\pm {
v}\}$ for all $ v\in R$.

The Coxeter group $G$ (or the f\/inite ref\/lection group) generated by
the root system $R$ is the subgroup of the orthogonal group $O(n)$
generated by $\{\sigma_{u}: u\in R\}$.

The positive subsystem $R_+$ is a subset of $R$ such that $R=R_+\cup
(-R_+)$, where $R_+$ and $-R_+$ are separated by some hyperplane
through the origin.

A multiplicity function
\begin{gather*}
\kappa : \ R\longrightarrow \mathbb C,
\\
\phantom{\kappa : {}}{} \  v \longmapsto \kappa_{v}
\end{gather*}
 is a $G$-invariant complex valued function def\/ined on $R$,
i.e.\ $\kappa_{v}=\kappa_{g{v}}$ for all $g\in G$.

 Fix a positive
subsystem $R_+$ of $R$ and denote
\[
\gamma=\gamma_\kappa:=\sum_{ v\in R_+}\kappa_{ v}.
\]

  Let $\mathcal{D}_j$ be the Dunkl operator associated to
the  group $G$ and to the multiplicity function $\kappa$, def\/ined by
 \begin{equation}\label{eq2.1}
 \mathcal{D}_jf(x)=\frac{\partial}{\partial
 x_j}f(x)+\sum_{\upsilon\in R_+}
 \kappa_\upsilon \frac{f(x)-f(\sigma_\upsilon x)}{\langle x,\upsilon\rangle}\upsilon_j.
 \end{equation}
We call $\Delta_h=\sum\limits_{j=1}^n {\mathcal{D}_j}^2$ the Dunkl
Laplacian. We  always assume that $ \kappa_\upsilon\ge 0. $

Let $d\sigma$ be the Lebesgue surface measure in the unit sphere and
$h_{k}(x)=\prod\limits_{\upsilon\in
R_{+}}|\langle\upsilon,x\rangle|^{\kappa_{\upsilon}}.$ Denote
$f^+=\frac{|f|+f}{2}$ and
\[
M_1 (r,f)=\int_{|y|=r}|f(y)|h_k ^2 (y)d\sigma(y).
\]

The mean value property holds for Dunkl harmonic functions $f$,
i.e.,
 \begin{equation}\label{eq2.2}
 \int_{|x|=1} f(x)h_{k}^{2}(x) d\sigma(x)=cf(0)\end{equation}
for some  constant $c>0$. This property for polynomials is implicit
in the orthogonality relation in~\cite{D}. Then one only needs a
limiting argument for arbitrary Dunkl harmonic functions. See also~\cite{MT, MY}.

Our main theorem is as follows.

\begin{theorem}\label{thm:main}
 Assume that $s\in\mathbb N\cup\{0\}$, $p\in\mathbb N$, and $s\geq2(p-1)$. Let
$f\in C^{2p}(\mathbb{R}^n)$ and $\Delta_{h}^p f=0$. Then $f$ is a
polynomial of degree $\le s$ if and only if
\[
\liminf_{r \rightarrow \infty} \frac{M_1(r,f^{+})}{r^{s+n-1+2 \gamma}} \in [0,+ \infty).
\]
Moreover, when  $s>2(p-1)$,  we have

 $(i)$ $f$ is a polynomial of
degree $<s$ if and only if
\[
\liminf_{r \rightarrow \infty} \frac{M_1(r,f^{+})}{r^{s+n-1+2 \gamma}}=0;
\]

$(ii)$ $f$ is a polynomial of degree $s$ if and only if
\[
\liminf_{r \rightarrow \infty} \frac{M_1(r,f^{+})}{r^{s+n-1+2 \gamma}} \in (0,+ \infty).
\]
\end{theorem}

As a direct corollary,  a Dunkl harmonic function on
$\mathbb{R}^{n}$ must be  constant if it is
 bounded below or above. Indeed, let $f$ be Dunkl harmonic in $\mathbb R^n$ and set
$s=0$.
 If $f\geq0$,  then $f^{+}=f$ and so that the mean value property
 shows that
\[
\frac{M_1(r,f^{+})}{r^{n-1+2 \gamma}}=\int_{|x|=1} f(x)h_{k}^{2}(x) d\sigma(x)=cf(0)\in [0,+ \infty).
\]
Therefore, Theorem~\ref{thm:main} shows that $f$ is a polynomial of
degree less than or equal to $s=0$, which means that $f$ must be a
constant. If $f$ is bounded above or below, multiplication by $-1$
if necessary makes it bounded below, adding a constant if necessary
makes it positive. Thus the known result for positive functions
shows that $f$ is a constant.

\section{Homogeneous expansions}

The non-Dunkl case of the following homogeneous expansions is
well-known (see \cite[Corollary~5.34]{ABR}). Its Dunkl version will
be helpful in the proof of Theorem \ref{thm:main}.

Denote by $\mathcal{ H}_m^n(h_\kappa^2)$ the space of Dunkl harmonic
polynomials of degree $m$ in $\mathbb R^n$.

\begin{lemma}\label{le:main}
If  $u$ is a  Dunkl harmonic function in $\mathbb R^n$, then there
exist $p_m\in \mathcal{ H}_m^n(h_\kappa^2)$, such that
\[
u(x)=\sum_{m=1}^{\infty} p_{m}(x), \qquad |x|<1,
\]
the series converging absolutely and uniformly on compact subsets of
the unit ball.
\end{lemma}

\begin{proof}
The formula can be verif\/ied as in the classical case of Corollary~5.34 in \cite{ABR} by using the formulae corresponding to the
classical case in~\cite{DX}. Indeed, notice that $u$ is harmonic on
the closed unit ball. Theorem~5.31 in~\cite{DX} and the preceding
statement  show that
\[
u(x)=\sum_{m=1}^\infty  c_h'\int_{|y|=1} u(y) P_m(h_\kappa^2; x, y)
h_\kappa^2(y) d\sigma(y),\qquad |x|<1,
\]
 where
\[
P_m(h_\kappa^2; x, y)=\sum_{0\le j \le n/2}
\frac{(\gamma_\kappa+\frac{n}{2})_m
2^{m-2j}}{(2-m-\gamma_\kappa-n/2)_jj!}|x|^{2j}|y|^{2k} K_{m-2j}(x,
y).
\] Let
\[
p_m(x)=c_h'\int_{|y|=1} u(y) P_m(h_\kappa^2; x, y)
h_\kappa^2(y) d\sigma(y),\qquad x\in\mathbb R^n.
\]
 Then
$p_m\in \mathcal{ H}_m^n(h_\kappa^2)$, since $P_m(h_\kappa^2; x, y)$
is the reproducing kernel of $\mathcal{ H}_m^n(h_\kappa^2)$ (see
\cite[p.~131, p.~189]{DX}).

By Proposition 4.6.2(ii) in~\cite{DX}, we have
\[
K_m(x, y)\le
\frac{1}{m!} |x|^m|y|^m,
\]
 so that
\[
|p_m(x)|\le C m^{2+\gamma+n/2} |x|^{m}\int_{|y|=1} |u(y)|
h_\kappa^2(y) d\sigma(y), \qquad \forall \, x\in\mathbb R^n,
\] and
thus the series $\sum_m p_m$ converges absolutely and uniformly to
$u$ on compact subsets of the unit ball.
\end{proof}

\section{Proof of the theorem}
\begin{proof} The necessity of Theorem \ref{thm:main} is clearly true. Now we prove  that $f$ is a polynomial of degree
$\leq s$ whenever
\begin{equation}\label{eq4.1}\liminf_{r\rightarrow\infty}\frac{M_1
(r,f^+)}{r^{s+n-1+2 \gamma}} \in [0,\infty).\end{equation} From the
Almansi decomposition theorem for Dunkl polyharmonic functions
\cite{R}, we have
\[
f(x)=\sum_{m=0}^{p-1}|x|^{2m}\varphi_m (x),
\]
 $\varphi_m$ being
Dunkl harmonic functions. Therefore, the mean value property in
(\ref{eq2.2}) shows
\begin{gather*}
 \int_{|y|=r}f(y)h_k^2 (y)d\sigma (y)  =
 \int_{|y|=r}  \sum_{m=0}^{p-1} |y|^{2m}
\varphi _m (y)h_{k}^2(y)d\sigma(y) \\
\phantom{\int_{|y|=r}f(y)h_k^2 (y)d\sigma (y)}{}
=\sum_{m=0}^{p-1}r^{2m} \int_{|y|=r}
\varphi _m (y)h_{k}^2
 (y)d\sigma(y)
 \\
\phantom{\int_{|y|=r}f(y)h_k^2 (y)d\sigma (y)}{} = \sum_{m=0}^{p-1}r^{2m+n-1+2\gamma} c\varphi_{m}(0)
= O(r^{2(p-1)+n-1+2\gamma}), \qquad
{r\rightarrow\infty}.
\end{gather*}
Thus
\[
\lim_{r\rightarrow\infty}\frac{1}{r^{j+n-1+2\gamma}}\displaystyle\int_{|y|=r}f(y)h_{k}^2(y)d\sigma(y)=
\begin{cases} c \, \varphi_{p-1}(0),\quad &j=2(p-1),\\
0, \quad \quad &j>2(p-1).\end{cases}
\]
Since $|f|=2f^+ -f$, we have
\begin{gather}
0 \leq   \liminf_{r\rightarrow\infty}\frac{M_1
(r,f)}{r^{s+n-1+2 \gamma}}   =2
\liminf_{r\rightarrow\infty}\frac{M_1 (r,f^{+})}{r^{s+n-1+2
\gamma}}- \lim_{r\rightarrow \infty}\frac{1}{r^{s+n-1+2
\gamma}}  \int_{|y|=r}f(y)h_{k}^2 (y)d\sigma (y)\nonumber\\
\phantom{0}{} =2  \liminf_{r\rightarrow\infty}\frac{M_1
(r,f^{+})}{r^{s+n-1+2 \gamma}}-C_1<+\infty, \label{eq4.30}
\end{gather}
where $C_1=c\varphi_{p-1}(0)$ for $s=2(p-1)$ and $C_1=0$ for
$s>2(p-1)$.

Since $\varphi_m$ are Dunkl harmonic  in $\mathbb R^n$, by Lemma
\ref{le:main}
 we can write
\[
\varphi_m=\sum_{j=1}^{\infty} g_{m,j},
\]
where $g_{m,j}$ are Dunkl harmonic homogeneous polynomials of degree
$j$ and the convergence of the series is uniform on compact subsets
of the unit ball.

We claim that $g_{m,j}=0$ when $2m+j>s$. With this claim, we have
\begin{equation}\label{eq4.2}
f(x)=\sum_{m=0}^{p-1}|x|^{2m}\varphi_{m}(x)=\sum_{m=0}^{p-1}|x|^{2m}\sum_{j=0}^{s-2m}g_{m,j},
\end{equation}
so that  $f(x)$ is a polynomial of degree no more than~$s$.

We now  prove  the claim by contradiction. Suppose there is a
$g_{m_{0},j_{0}}$ not identical $0$ and $2m_{0}+j_0 >s $, then
\begin{gather}
  \liminf_{r\rightarrow\infty}
\frac{1}{r^{s+j_0+n-1+2\gamma}}\displaystyle
\left|\int_{|y|=r}f(y)g_{m_{0},j_{0}}(y)h_{k}^2(y) d\sigma (y)\right|\nonumber\\
\qquad{}\leq \liminf_{r\rightarrow\infty}
 \frac{1}{r^{s+j_0+n-1+2\gamma}}
\max_{|y|=r}|g_{m_{0},j_{0}}(y)|
\int_{|y|=r}|f(y)|h_{k}^2 (y)d\sigma (y)\nonumber\\
 \qquad{}= \max_{|y|=1}|g_{m_{0},j_{0}}(y)|\displaystyle
\liminf_{r\rightarrow\infty} \frac{M_1 (r,f)}{r^{s+n-1+2\gamma }}\in
[0,\infty). \label{eq4.3}
\end{gather}

On the other side, notice that $\{g_{m,j}\}_{j=0}^\infty$ are
orthogonal in $L^2(S_{n-1}, h_k^2 d\sigma)$ (see~\cite{DX}), where
$S_{n-1}$ is the unit sphere in $\mathbb R^n$, we have
\begin{gather*}
 \int_{|y|=r}f(y)g_{m_{0},j_{0}}(y)h_{k}^2(y) d\sigma
(y)=\int_{|y|=r}\displaystyle\sum_{m=0}^{p-1}|y|^{2m}\displaystyle\sum_{j=0}^{\infty}g_{m,j}(y)g_{m_0
,j_0}(y)h_k ^2 (y)d\sigma(y)
\\
 \qquad{} = \int_{|y|=r} \sum_{m=0}^{p-1}|y|^{2m}g_{m,j_0}(y)g_{m_0
,j_0}(y)h_k ^2 (y)d\sigma(y).
\end{gather*}
On considering item $m=m_0$, the above integral is of the form
$r^{2\gamma}P(r)$ with $P(r)$ a polynomial of $r$ with degree at
least $2m_0 +2j_0 +n-1$, which is strictly larger than $s+j_0 +n-1$.
This contradicts the fact that
\begin{gather}
\displaystyle\int_{|y|=r}f(y)g_{m_{0},j_{0}}(y)h_{k}^2 (y) d\sigma
(y) = \int_{|y|=r} (2f^+(y)-|f(y)|)
g_{m_{0},j_{0}}(y)h_{k}^2 (y) d\sigma (y)\nonumber\\ \qquad =
O(r^{s+j_0+n-1+2\gamma}). \label{eq4.4}
\end{gather}
The last step used (\ref{eq4.1}) and (\ref{eq4.3}).

We have proved that   $f$ is a polynomial of degree $\leq s$ if and
only if (\ref{eq4.1}) holds.  With this result, we know that the
statements $(i)$ and $(ii)$ in Theorem \ref{thm:main} are equivalent.

It remains to prove the suf\/ficiency of Theorem \ref{thm:main}$(i)$ in
the case $s>2p-2$.

Let $s>2p-2$ and assume that $f$ satisfies
\[
\liminf_{r
\rightarrow \infty} \frac{M_1(r,f^{+})}{r^{s+n-1+2 \gamma}}=0.
\]
 In
this case, from the proof of  (\ref{eq4.30}), we f\/ind
\[
\liminf_{r\rightarrow\infty}\frac{M_1 (r,f)}{r^{s+n-1+2
\gamma}}=0,
\]
 so that the same reasoning as in (\ref{eq4.4}) shows
that
\[
\int_{|y|=r}f(y)g_{m,j}(y)h_{k}^2 d\sigma (y)=o(r^{s+n-1+j +2
\gamma}),\qquad r\to\infty.
\] We now claim that $g_{m,j}=0$ if
$2m+j\geq s.$

Indeed, if not, then there exists $g_{m_0 ,j_0}$ $\not\equiv 0$ such
that $2m_0 +j_0 \geq s$. However,
\[
\int_{|y|=r}f(y)g_{m_{0},j_{0}}(y)h_{k}^2 d\sigma (y)
=\int_{|y|=r}\sum_{m=0}^{p-1}|y|^{2m}\sum_{j=0}^{\infty}g_{m,j}(y)g_{m_0
,j_0}(y)h_k ^2 (y)d\sigma(y),
\]
and the above integral is of the form $r^{2\gamma}Q(r)$, where
$Q(r)$ is  a polynomial of $r$ with degree no less than $2m_0 +2j_0
+n-1\ge s +j_0 +n-1$. Since $r^{2\gamma}Q(r)$ is not an
$o(r^{s+j_0+n-1 +2 \gamma})$, we arrive at a contradiction.

From the claim and (\ref{eq4.2}), we see that $f$ is a polynomials
of degree less than $s$.
\end{proof}

\subsection{Further remarks}

We have shown that our main results hold for any non-negative
multiplicity function $\kappa$. With the same approach, we can
extend our main result to the case when
\[
\kappa\in K^{\rm reg} \qquad {\mathrm{and}}\qquad  \Re \mathfrak{e} \,\kappa>-\frac{n}{2},
\]
where $K^{\rm reg}$ is the regular parameter set (see \cite{DJO, Ro}).
This is because in the proof we applied the Dunkl intertwining
operator in Lemma~\ref{le:main} and the Almansi decomposition. The
existence of the intertwining operator depends on the restriction
that $\kappa\in K^{\rm reg}$,  while the Almansi decomposition holds for
any $\Re \mathfrak{e}\, \kappa>-\frac{n}{2}$ (see \cite{R, R1}). We
mention that in the rank-one case, i.e., \mbox{$\dim
\mathrm{span}_{\mathbb R}R=1$}, we have
$
K^{\rm reg}=\mathbb C\setminus\{-1/2-m, \   m=0, 1, 2,  \ldots\}.
$
This means  that in the rank-one case, Theorem~\ref{thm:main} holds
for
$
\Re \mathfrak{e} \,\kappa>-\frac{1}{2}.
$

\subsection*{Acknowledgements}

The authors would like to thank  the  referees for their useful
comments.  The research is supported by the Unidade de
Investiga\c c\~ao ``Matem\'atica e Aplica\c c\~oes'' of University
of Aveiro, and by the NNSF  of China
 (No. 10771201), NCET-05-0539.

\pdfbookmark[1]{References}{ref}

\LastPageEnding

\end{document}